\documentclass[a4paper]{article}

\newcommand{\pic}[2]{\BoxedEPSF{#1 scaled #2}}
\def\meridianmap{\pic{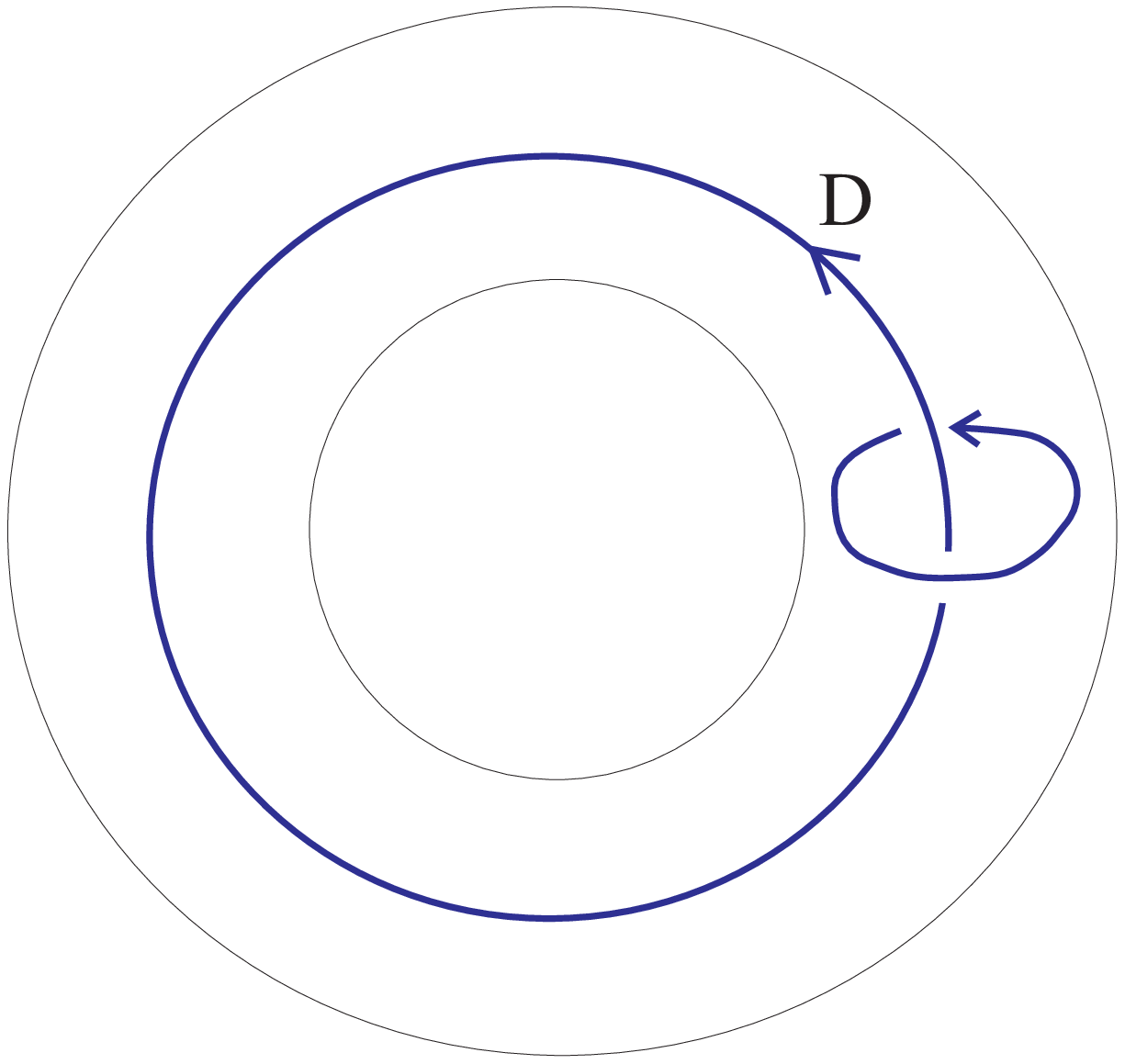}{300}}
\def\npbox{\pic{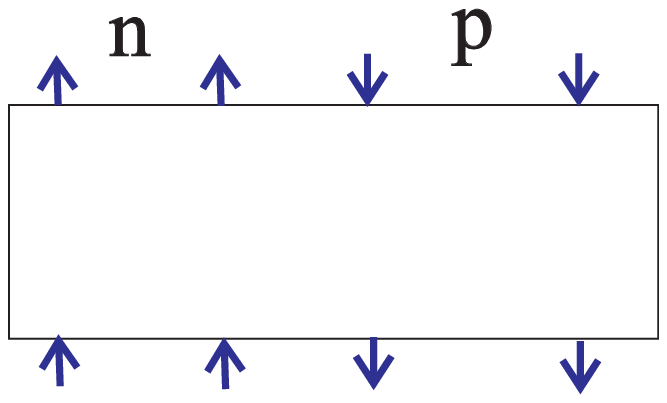}{400}}
\def\onetangle{\pic{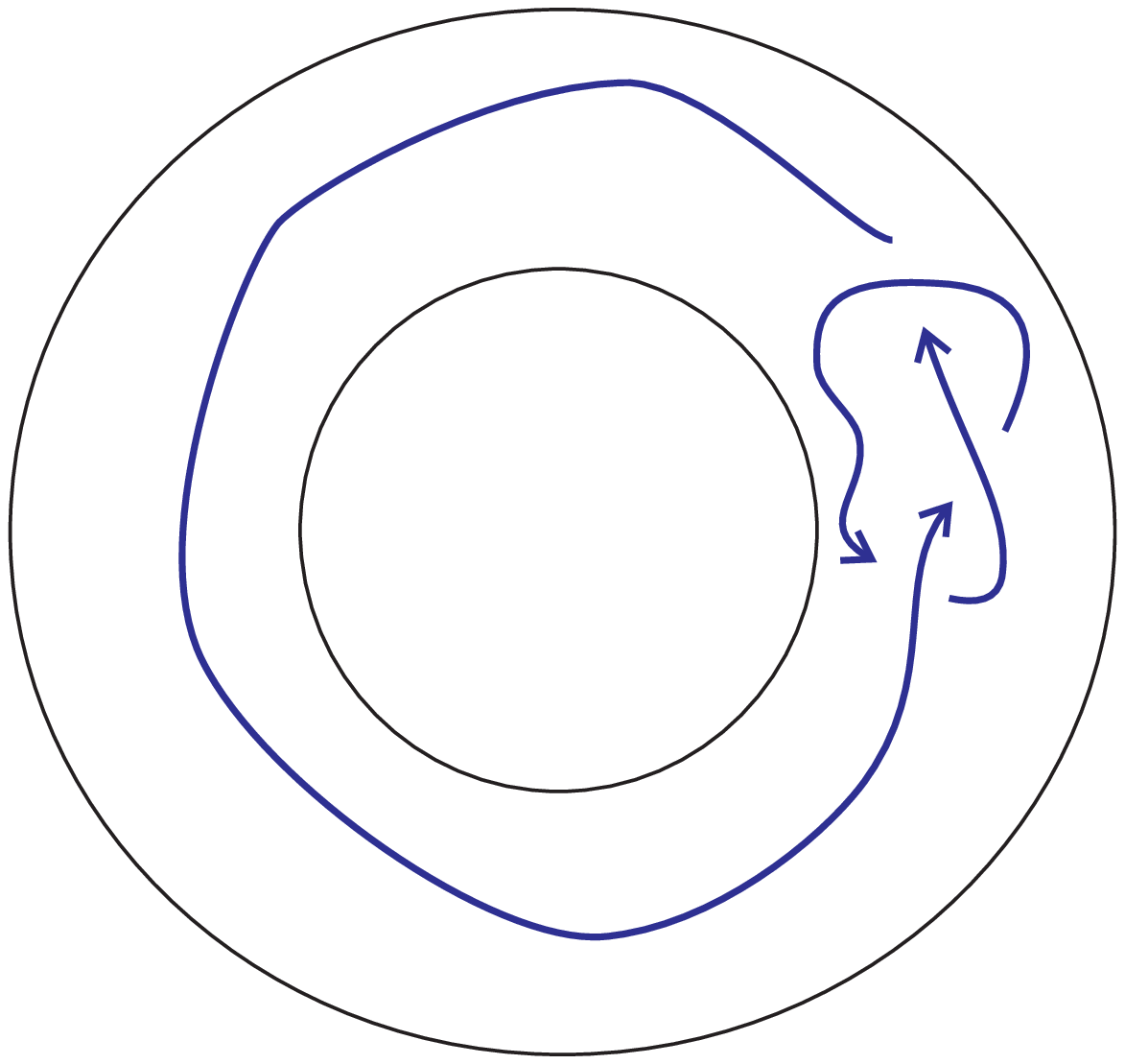}{300}}
\def\onetangleQ{\pic{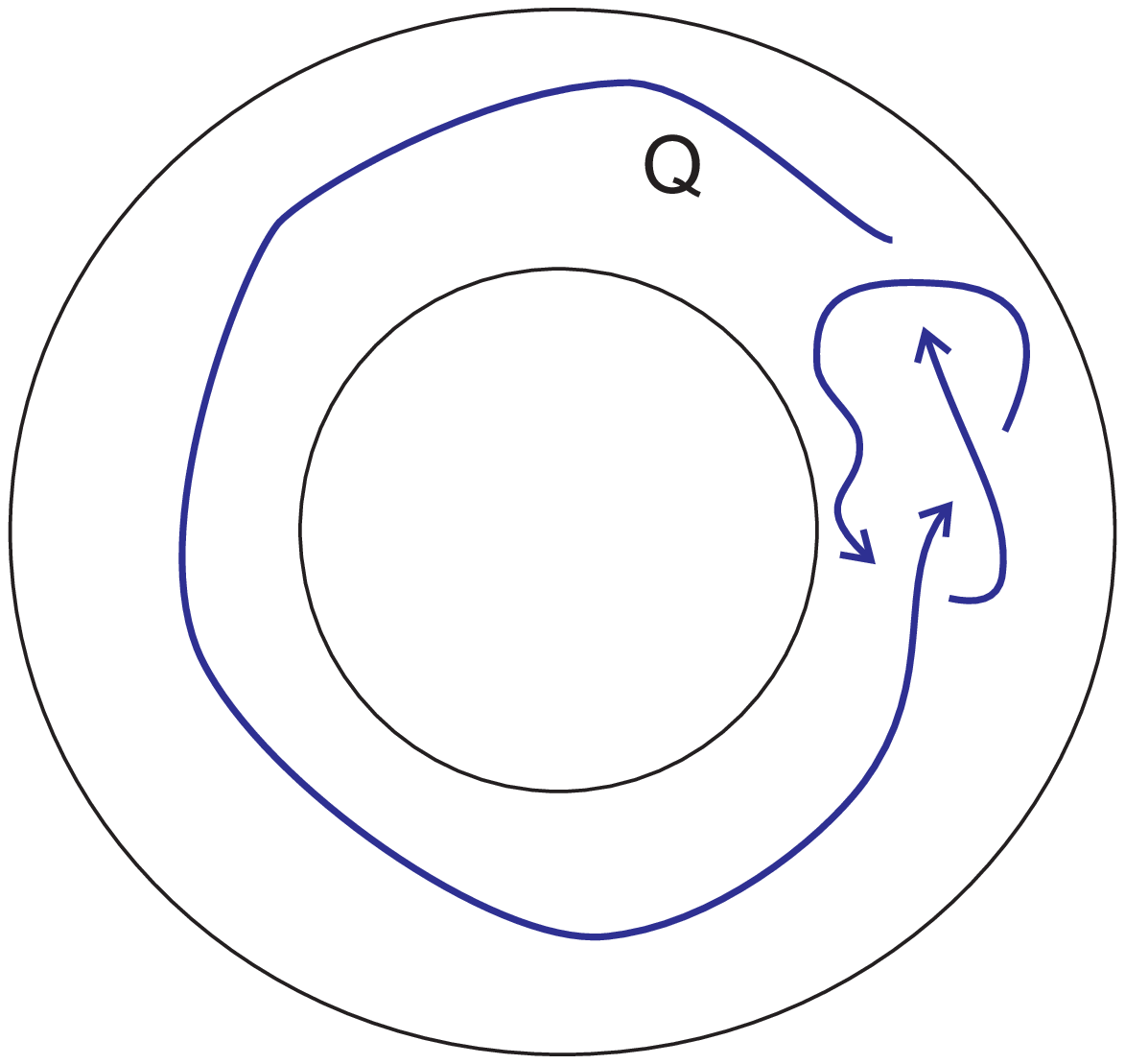}{300}}

\def\oneonebox{\pic{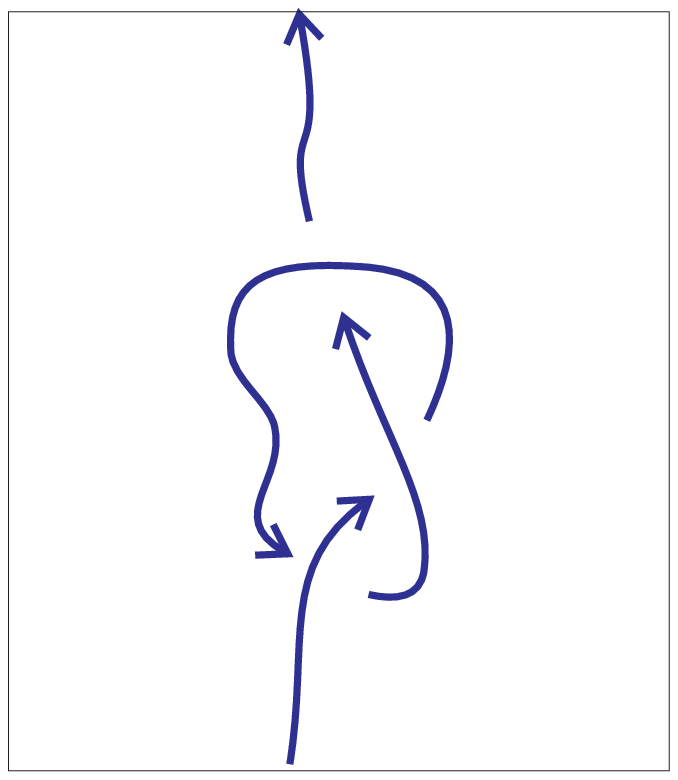}{300}}
\def\twooneparallel{\pic{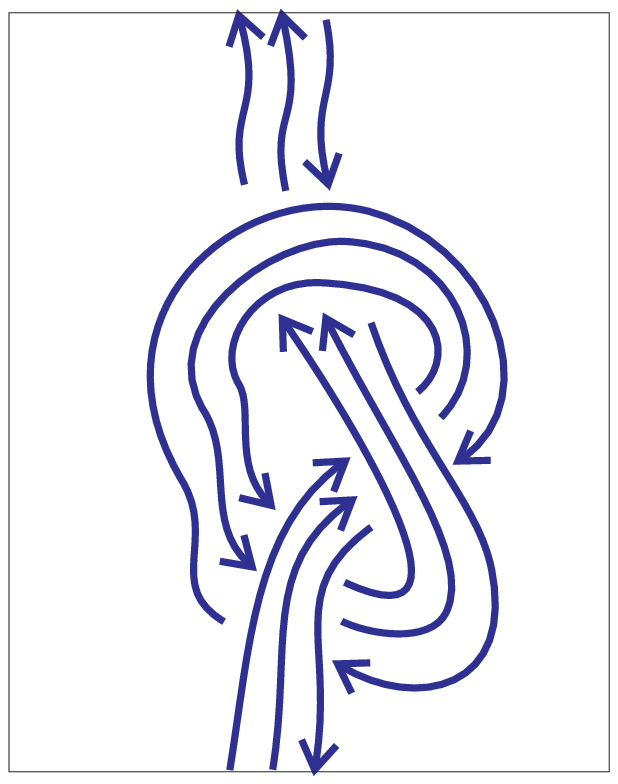}{400}}

\def\Xor{\pic{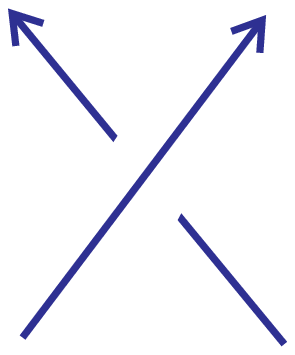} {250}}
\def\Yor{\pic{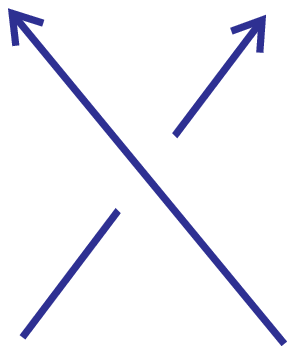} {250}}
\def\Ior{\pic{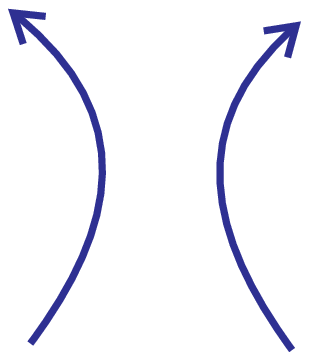} {250}}

\def\unknot{\pic{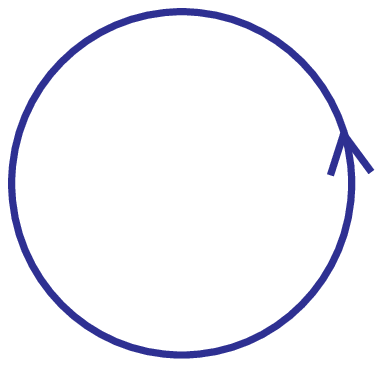} {150}}

\def\Idor{\pic{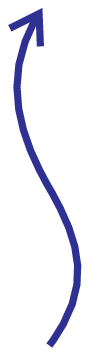} {250}}
\def\Rcurlor{\pic{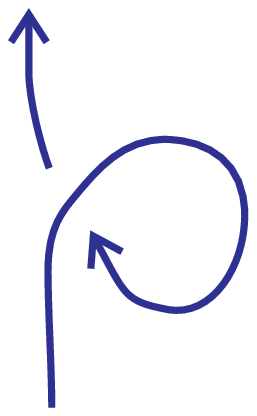} {250}}

\def\CC{\mathcal {C}}

\newcommand{\bibname}{}

\usepackage[usenames]{color}

\usepackage{multicol}

\input{dvipsnam.def}

\include{symboldefs}

\usepackage[hideboxes]{boxedeps}
\usepackage{amsmath}
\usepackage{amssymb}

\newcommand{\bc}{\begin{center}}
\newcommand{\ec}{\end{center}}
\newcommand{\gl}{\lambda}
\newcommand{\qlm}{Q_{\gl,\mu}}
\newcommand{\ds}{\displaystyle}

\newtheorem{theorem}{Theorem}
\newtheorem{lemma}{Lemma}
\newtheorem{conjecture}{Conjecture}

\newenvironment{remark}{\par\smallskip%
\noindent\textbf{Remark.}\  }%
{\par\smallskip}

\newenvironment{proof}[1][{}]{\par\smallskip%
\noindent\textit{Proof #1: }\  }
{\hfill$\Box$\par\smallskip}

\begin{document}

\title{
  Integrality of Homfly $(1,1)$-tangle invariants}

\author{H. R. Morton }

\maketitle
\begin{abstract}
Given an invariant $J(K)$ of a knot $K$, the corresponding $(1,1)$-tangle invariant $J'(K)={J(K)}/{J(U)}$ is defined as the quotient  of $J(K)$ by its value $J(U)$ on the unknot $U$.
We prove here that $J'$ is always an integer 2-variable Laurent polynomial when $J$ is the Homfly satellite invariant determined by decorating  $K$ with any eigenvector of the meridian map in the Homfly skein of the annulus. Specialisation of the 2-variable polynomials for suitable choices of eigenvector shows that the $(1,1)$-tangle irreducible quantum $sl(N)$ invariants of $K$ are integer 1-variable Laurent polynomials.
\end{abstract}

\section*{Introduction} 

Decorating a framed knot $K$ with a pattern $Q$ (a diagram in the standard annulus) determines a satellite $K*Q$ of $K$, whose Homfly polynomial is a 2-variable Laurent polynomial $P({K*Q})\in{\bf Z}[v^{\pm1},z^{\pm1}]$.  For each fixed $Q$ this gives a 2-variable invariant of the knot $K$. We admit linear combinations of patterns, regarded as elements of the Homfly skein of the annulus, in place of single diagrams $Q$, and extend our coefficients to the ring $\Lambda$ of  Laurent polynomials ${\bf Z}[v^{\pm1},s^{\pm1}]$ with denominators $s^r-s^{-r}, r\ge1$, taking $z=s-s^{-1}$, to provide an invariant $J(K)=P({K*Q})\in\Lambda$ for any $\Lambda$-linear combination $Q$ of patterns.

For each partition $\lambda$ of $n$ and each $N$, the quantum $sl(N)_q$ invariant of $K$ when colored by the irreducible module corresponding to  $\gl$ is an integral Laurent polynomial in $s$, with $q=s^2$.  It has been known for some time, \cite{Wenzl,AistonMorton, Kawagoe, Lukac}, how to choose a decoration $Q_\lambda$  so that the 2-variable Homfly invariant $P({K*Q_\gl})$ gives all these 1-variable invariants for different values of $N$ by substituting $v=s^N$.  The invariant $P({K*Q_\gl})$ typically involves denominators $s^r-s^{-r}$ with $r$ up to the maximum hook-length of the partition $\gl$.

In \cite{Le} Thang Le showed that the $(1,1)$-tangle invariant $J'_K(V_\gl)$ of a framed knot $K$ when colored by an irreducible module $V_\gl$ over any quantum group is an integer Laurent polynomial in the quantum parameter $q$.  In this case the `quantum dimension' of  $V_\gl$, which is  $J_U(V_\gl)$, is itself in ${\bf Z}[q^{\pm1}]$ and hence so is the invariant $J_K(V_\gl)$.

Consequently the denominators in the 2-variable invariant $P({K*Q_\gl})$ will be cancelled by terms in the numerator when $v$ is replaced by $s^N$ for any $N$.
Each of the resulting 1-variable Laurent polynomial invariants of $K$ is then divisible by the value of the invariant for the unknot. Many constructions of manifold invariants based on quantum invariants involve substitution of a root of unity for the variable $s$; the $(1,1)$-tangle invariant gives a reliable means of retaining information at values of $s$ for which the quantum dimension of the coloring module is zero.

The purpose of this paper is to show that the integrality of these $sl(N)_q$  $(1,1)$-tangle invariants of $K$ can already be seen at the 2-variable level. We show further that if $J(K)=P(K*Q)$ where $Q$  is \emph{any} eigenvector of the meridian map on the Homfly skein of the annulus then the $(1,1)$-tangle invariant $J'(K)$ lies in  ${\bf Z}[v^{\pm1},s^{\pm1}]$. Such eigenvectors $Q$  include the elements $Q_\gl$ mentioned already, as well as a wider family $\qlm$, \cite{HadjiMorton}, depending on two partitions $\gl$ of $n$ and $\mu$ of $p$. These give a single 2-variable invariant which packages together for different $N$ the quantum invariants coming from the irreducible submodule of the tensor product of $n$ copies of the fundamental $sl(N)_q$ module and $p$ copies of its dual determined by the partitions $\gl$ and $\mu$. The individual 1-variable invariants are recovered from $P(K*\qlm)$ in the form of a single 2-variable integral invariant $J'(K)=a_K(\gl,\mu)$ which yields each $sl(N)_q$ invariant by setting $v=s^N$. In the simplest case where $n=p=1$ the modules are the adjoint representations of $sl(N)_q$, and the 2-variable invariant is closely related to the Homfly polynomial of the reverse parallel of the knot.

The eigenvectors $\qlm$ of the meridian map in the Homfly skein of the annulus are described explicitly in \cite{HadjiMorton}, where further details of their properties can be found. The main result here is the following integrality theorem for the 2-variable $(1,1)$-tangle invariants $\ds a_K(\gl,\mu)=\frac{P(K*\qlm)}{P(U*\qlm)}$ of a framed knot $K$ coming from $J(K)=P(K*\qlm)$.

\begin{theorem} Let $K$ be a framed knot and
let $Q$ be any eigenvector of the meridian map. Then the $(1,1)$-tangle invariant $a_K=P(K*Q)/P(U*Q)$ is a 2-variable integer Laurent polynomial $a_K\in{\bf Z}[v^{\pm1},s^{\pm1}]$.
\end{theorem}
\setcounter{theorem}{0}

As a corollary the Homfly polynomial $P(K*\qlm)$ of the satellite $K*\qlm$ will always factorise  as $P(K*\qlm)=a_K(\gl,\mu)P(U*\qlm)$ with $a_K(\gl,\mu)\in {\bf Z}[v^{\pm1},s^{\pm1}]$.

The proof depends on controlling the powers of $z^{-1}$ in a skein resolution of a single diagram in a surface in terms of the number of null-homotopic closed components of the diagram. Calculations  in which braids interact with an element of the Hecke algebra which closes to give $Q_\gl$, based on \cite{AistonMorton}, are then combined with relations from \cite{MortonHadji} between $\qlm$, $Q_\gl$ and $Q_\mu$ to complete the argument.

\section{Homfly skeins and resolutions}

\subsection*{The general setting}

Homfly skein theory applies to a surface $F$ with some distinguished input and output boundary points.

The (linear) {skein} of $F$ is defined as linear combinations of diagrams in $F$, up to Reidemeister moves II and III,  modulo the skein relations
\bc\begin{enumerate}
\item \qquad{$\Xor\  -\ \Yor \qquad =\qquad{(s-s^{-1})}\quad\ \Ior \ ,$}
\item \qquad
{$ \Rcurlor \qquad=\qquad {v^{-1}}\quad \Idor\ . $}
\end{enumerate}
\ec

The coefficient ring $\Lambda$ is  taken as $Z[v^{\pm1},s^{\pm1}]$, with denominators $\{r\}=s^r-s^{-r}, r\ge1$.

 Application of the first relation to the crossing in the second relation gives the relation  $(v^{-1}-v)\ \Idor\quad=\quad z\ \unknot\ \Idor\ $. This can be used to remove a null-homotopic curve $\unknot$ without crossings from a diagram at the expense of introducing  $z^{-1}$ in the coefficients. 

\subsection*{Examples}
The skein of the plane  is spanned by a single element, \unknot.
Any link $L$ represents $P(L)\ \unknot$ where $P(L)\in\Lambda$ is its Homfly polynomial.

When $F$ is a rectangle with $n$ outputs and $p$ inputs at the top, matched at the bottom as shown

\bc \npbox
\ec 
the diagrams are called $(n,p)$-tangles, although to be consistent with the  terminology for $(1,1)$-tangles their name might be expanded to $([n,p],[n,p])$-tangles.

The resulting skein, $H_{n,p}$, has finite dimension $(n+p)!$, and is  an {algebra} over $\Lambda$, where the product is induced by placing one tangle above another.

\subsection*{Resolutions}

A \emph{resolution tree} for a diagram $D$ in $F$ is a directed tree of diagrams in $F$, with initial vertex $D$, having either one or two edges leaving each internal vertex. Two edges lead to the diagrams where one crossing in the current diagram is either switched or smoothed. A single edge performs a Reidemeister move of type I on the current diagram or removes  a null-homotopic closed curve without crossings.

The following general integral resolution lemma controls the use of negative powers of $z$, and will shortly be applied in  $H_{n,p}$.
Write $ k(D)$ for the number of \emph{null-homotopic} closed curves in a diagram $D$. 

\begin{lemma} Let $D$ be a diagram in a surface $F$ having a resolution tree with diagrams $\{D_i:i\in I\}$ at its end vertices. Then $D$ can be written in the skein of $F$ as a $\Lambda$-linear combination of $\{D_i\}$ in the form
\[z^{k(D)}D=\sum_{i\in I}c_i z^{k(D_i)}D_i,\] where $c_i\in{\bf Z}[v^{\pm1},z]$.
\end{lemma}

\begin{proof} By induction on the number of edges of the resolution tree.  

\begin{itemize}
\item{\bf 1.}\quad
If two edges leave the vertex $D$ then the resolution has switched or smoothed a crossing of sign $\pm1$ in $D$,  resulting in diagrams $D_\mp$ and $D_0$ which satisfy $D=D_\mp \pm zD_0$. Now $k(D_\mp)=k(D)$ while $k(D_0)\le k(D)+1$. Then
\[z^{k(D)}D=z^{k(D_\mp)}D_\mp  \pm z^a z^{k(D_0)}D_0,\] with $a\ge 1$, and the resolution subtrees for $D_\mp$ and $D_0$ allow the right hand side to be expanded in terms of the end vertices $D_i$ by induction. The coefficients $c_i$ are either unchanged or multiplied by $\mp z^a, a\ge0$.

\item{\bf 2.}\quad
If a single edge leaving $D$ comes from a Reidemeister  type I move then the result is immediate. If the edge corresponds to the removal of a null-homotopic closed curve without crossings, leading to a diagram $D'$, then $k(D')=k(D)-1$, while $zD=(v^{-1}-v)D'$ in the skein. Then \[z^{k(D)}D=(v^{-1}-v)z^{k(D')}D',\] and again induction gives the required expansion, using the subtree for $D'$ whose coefficients $c_i$ are multiplied by $(v^{-1}-v)$.
\end{itemize}
The induction starts trivially for a resolution tree with $0$ edges.
\end{proof}

\subsection*{Resolutions in $H_{n,p}$}

A framed knot $K$ can be represented as a $(1,1)$-tangle $T(K)$ by a single knotted arc, such as \[\oneonebox\ .\] The $(n,p)$-parallel of this, $T_{n,p}(K)$, in the skein $H_{n,p}$ is constructed by drawing $n+p$ parallel strands to the arc $T(K)$, with $n$ oriented in one sense and $p$ in the other, illustrated here with $n=2,p=1$,
\[\twooneparallel\ .\]
 Standard procedures allow its resolution into $(n+p)!$ \emph{totally descending} tangles without closed components; these are tangles in which every crossing is first met as an overcrossing when the arcs are traversed in order. The ordering of the arcs in each of these tangles can be chosen by ordering their initial points counterclockwise around the boundary, starting from the bottom left corner. As a corollary of the integrality lemma above, $T_{n,p}(K)$ can be written as a linear combination  of these tangles with all coefficients in ${\bf Z}[v^{\pm1},z]$.

In the case $p=0$ such tangles are the `positive permutation braids', $\{b_\pi;\pi\in S_n\}$, with strings oriented from bottom to top, while when $n=0$ they are again positive permutation braids $\{b^*_\rho;\rho\in S_p\}$, with string orientation from top to bottom. In general each tangle is determined up to isotopy by knowing which input and output points are connected by its arcs. 

 For each tangle we may count the number $k$ of its arcs which connect input and output points at the bottom. Then $0\le k\le \min(n,p)$. We can write $T_{n,p}(K)$ in the skein $H_{n,p}$ as $T_{n,p}(K)=T^{(0)}_{n,p}(K)+T^{(0)}_{n,p}(K)$ where $T^{(0)}_{n,p}(K)$ is a combination of tangles with $k=0$ and $T^{(1)}_{n,p}(K)$ is a combination of tangles with $k\ge1$. Tangles with $k=0$ have the form $b_\pi\otimes b^*_\rho$ for some $\pi\in S_n$ and $\rho\in S_p$, where $\otimes$ denotes juxtaposition of tangles side by side. We then have
\[T^{(0)}_{n,p}(K)=\sum_{\pi\in S_n,\rho\in S_p}c_{\pi,\rho}(K)(b_\pi\otimes b^*_\rho),\] with all coefficients $c_{\pi,\rho}(K)$ in ${\bf Z}[v^{\pm1},z]$.

The subspace $H^{(1)}_{n,p}$ of the algebra $H_{n,p}$ spanned by the totally descending tangles with $k\ge 1$ forms a $2$-sided ideal, and indeed is one of a chain of ideals $H^{(l)}_{n,p}$, spanned by the tangles with $k\ge l$, which are discussed further in \cite{MortonHadji}. The closure map, induced by  taking an $(n,p)$-tangle to its closure in the annulus, carries the skein $H_{n,p}$ to a subspace $\CC_{n,p}$ of the skein $\CC$ of the annulus.  The image of  $H^{(1)}_{n,p}$ under this map can readily be seen to lie in $\CC_{n-1,p-1}\subset\CC_{n,p}$. In much of what follows we can work modulo $\CC_{n-1,p-1}$, so that the element $T^{(1)}_{n,p}(K)$ will not figure largely in the calculations.

\section{The meridian map}

The skein of the annulus $\CC$ has been studied extensively, starting with work of Turaev \cite{Turaev}. It forms a commutative algebra over $\Lambda$, with the product induced by placing two diagrams in concentric annuli.  The \emph{meridian map} $\varphi:\CC\to\CC$ is induced by including a single meridian curve around a diagram $D$  in the thickened annulus to give the diagram \[\varphi(D)\quad=\quad\meridianmap\ .\]

\subsection*{Satellites}
Diagrams in the annulus are sometimes known as \emph{patterns} when they are used in the construction of satellites of a framed knot. Starting with a framed knot $K$ and a pattern $Q$, the satellite $K*Q$ is formed by replacing the framing annulus around $K$ with the annulus containing $Q$. This operation, known as \emph{decorating $K$ by $Q$},  induces a linear map at the skein level, so that the Homfly polynomial $P(K*Q)$ depends only on $ Q$ as an element of the skein $\CC$. If $K$ is drawn in the annulus as the closure of a  $(1,1)$-tangle 
\[\onetangle\] then decorating it by $Q$ gives a diagram of $K*Q$ in the annulus, \[\onetangleQ\] and induces a linear map $T_K:\CC\to \CC$. If $Q$ is an eigenvector of $T_K$ with eigenvalue $a_K$ then $K*Q=T_K(Q)=a_K Q=a_K U*Q$ where $U$ is the unknot with framing $0$. Taking the Homfly polynomial then gives $a_K=P(K*Q)/P(U*Q)$ as the $(1,1)$-tangle invariant $J'(K)$ coming from $J(K)=P(K*Q)$.

\subsection*{Eigenvectors}
The subspaces $\CC_{n,p}\subset \CC$ are invariant under the meridian map $\varphi$, and under $T_K$. A basis $Q_\gl$ of $\CC_{n,0}$ consisting of eigenvectors of $\varphi$ determined by partitions $\gl$ of $n$ has been described in \cite{AistonMorton}. The element $Q_\gl$ is constructed there as the closure of an idempotent $e_\gl$ in the skein $H_{n,0}$, which is isomorphic to the Hecke algebra $H_n$ of type A.  More recent constructions of Kawagoe and Lukac, \cite{Kawagoe,Lukac}, following the interpretation of $\CC_{n,0}$ as symmetric functions of degree $n$ in $N$ variables, show that the counterpart of the Schur function $s_\gl$ is also an eigenvector of the meridian map which can be identified with $Q_\gl$. The existence of a basis for the whole space $\CC$ consisting of eigenvectors of $\varphi$ with \emph{distinct} eigenvalues, indexed by pairs $\gl,\mu$ of partitions, is established in \cite{MortonHadji}, and explicit formulae for the eigenvectors $\qlm$ are given in \cite{HadjiMorton}.  Any eigenvector $Q$ of the meridian map is then a multiple of $\qlm$ for some partitions $\gl,\mu$.

\subsection*{Integrality}
We are now in a position to establish the main integrality result.
\begin{theorem} Let $K$ be a framed knot and
let $Q$ be any eigenvector of the meridian map. Then the $(1,1)$-tangle invariant $a_K=P(K*Q)/P(U*Q)$ is a 2-variable integer Laurent polynomial $a_K\in{\bf Z}[v^{\pm1},s^{\pm1}]$.
\end{theorem}

\begin{proof}
It is readily noted, \cite{HadjiMorton}, that the map $T_K$ commutes with the meridian map $\varphi$. Since the eigenvalues of $\varphi$ are distinct then any eigenvector $Q$ of $\varphi$ is also an eigenvector of $T_K$. The $(1,1)$-tangle invariant $J'(K)$ coming from the satellite invariant $J(K)=P(K*\qlm)$ is then the eigenvalue $a(\gl,\mu)$ of $T_K$ for its eigenvector $\qlm$. 
The integrality of $a(\gl,\mu)$ will now be established using features of $Q_\gl$ and $\qlm$ from \cite{AistonMorton} and \cite{MortonHadji}.

Turning the annulus over induces a symmetry $*$ in $\CC$ which carries an element $Q$ to $Q^*$. If $Q\in\CC_{n,p}$ then $Q^*\in\CC_{p,n}$. Thus if $\gl$ is a partition of $n$ and $\mu$ is a partition of $p$ we have $Q_\gl\in\CC_{n,0}$ and $Q_\mu^*\in\CC_{0,p}$ and their product $Q_\gl Q^*_\mu$ lies in $\CC_{n,p}$.  

In \cite{MortonHadji} it is shown that $\qlm=Q_\gl Q^*_\mu +W$ where $W\in\CC_{n-1,p-1}$.  Now $T_K(\qlm)=a(\gl,\mu)\qlm$ so $T_K(Q_\gl Q^*_\mu)=a(\gl,\mu)Q_\gl Q^*_\mu+V$ where $V\in\CC_{n-1,p-1}$.

The idempotent $e_\gl$   in \cite{AistonMorton}, whose closure is $Q_\gl$, can be factorised, following lemma 11 there, as $e_\gl=e^{(a)}_\gl e^{(b)}_\gl$ so that $e_\gl^{(a)}\beta e_\gl^{(b)}=k(\beta,\gl)e_\gl$ with $k(\beta,\gl)\in{\bf Z}[s^{\pm1}]$, for any $n$-braid $\beta$.
It follows that the closure of $e_\gl \gamma$, which is also the closure of $e_\gl \gamma e_\gl$, can be written as $c(\gamma,\gl) Q_\gl$, with $c(\gamma,\gl)\in{\bf Z}[s^{\pm1}]$, for any $n$-braid $\gamma$.

We can express $T_K(Q_\gl Q^*_\mu)$ as the closure of the element $(e_\gl\otimes e^*_\mu )T_{n,p}(K)$ in $H_{n,p}$.
Now \[(e_\gl\otimes e^*_\mu )T_{n,p}(K)=(e_\gl\otimes e^*_\mu )T^{(0)}_{n,p}(K), \bmod H^{(1)}_{n,p}.\]  
The closure of \[(e_\gl\otimes e^*_\mu )T^{(0)}_{n,p}(K)=\sum_{\pi\in S_n,\rho\in S_p}c_{\pi,\rho}(e_\gl b_\pi\otimes e^*_\mu b^*_\rho)\] is a scalar multiple  $A(\gl,\mu)Q_\gl Q^*_\mu$,  where \[A(\gl,\mu)=\sum_{\pi\in S_n,\rho\in S_p}c_{\pi,\rho}(K)c(b_\pi,\gl)c(b_\rho,\mu) \in {\bf Z}[v^{\pm1},s^{\pm1}]. \] 

Then $T_K(Q_\gl Q^*_\mu)=A(\gl,\mu)Q_\gl Q^*_\mu$ modulo $C_{n-1,p-1}$. Hence $A(\gl,\mu)=a(\gl,\mu)$ is the $(1,1)$-tangle invariant $P(K*\qlm)/P(U*\qlm)$, which is  a 2-variable integer Laurent polynomial in ${\bf Z}[v^{\pm1},s^{\pm1}]$, as claimed.
\end{proof}

\section{Some relations}

The $(1,1)$-tangle invariants $a(\gl,\mu)$ of $K$ are not all independent.

Firstly there are some symmetries.

\begin{itemize}
\item
By reversing orientation of all strings we get  $a(\mu,\lambda)=a(\lambda,\mu)$.

\item
Replacing $\lambda$ and $\mu$ by their {conjugate partitions} switches $s$ for $-s^{-1}$ in $a(\lambda,\mu)$.
\end{itemize}

Secondly the \emph{1-variable} invariant $\displaystyle a(\lambda,\mu)|_{{v=s^N}}$ agrees with $\displaystyle a(\nu)|_{{v=s^N}}$ for some explicit $\nu$ depending on $N,\lambda,\mu$, and corresponds to an {irreducible quantum $sl(N)$} invariant. Details of the appropriate partition $\nu$ can be found in \cite{HadjiMorton}.

An explicit determinantal construction for $\qlm$ is given in \cite{HadjiMorton} in terms of the elements $h_n=Q_{\lambda,\mu}$ where $p=0$ and $\lambda$ has a single part, and $h^*_p$ with the reverse orientation, where $n=0$ and $\mu$ has a single part. These elements generate $\CC$ freely as an algebra.

The general construction of $\qlm$ in \cite{HadjiMorton} can be  illustrated by the case when $\lambda$ has parts ${2,2,1}$ and $\mu$ has parts ${3,2}$.

Take a matrix with diagonal entries as shown, corresponding to the parts of $\lambda$ and $\mu$.

\[\begin{pmatrix} { h^*_2}&&&&\\
&{ h^*_3}&&&\\
&&{ h_2}&&\\
&&&{ h_2}&\\
&&&&{ h_1}
\end{pmatrix}\]

Complete the rows by shifting indices upwards for the parts of $\lambda$,
and downwards for the parts of $\mu$, to get
\[M=\begin{pmatrix} { h^*_2}&h^*_1&1&0&0\\
h^*_4&{ h^*_3}&h^*_2&h^*_1&1\\
1&h_1&{ h_2}&h_3&h_4\\
0&1&h_1&{ h_2}&h_3\\
0&0&0&1&{ h_1}
\end{pmatrix}\]

Then $Q_{\lambda,\mu}=\det M$.

\begin{remark}
The subalgebra of $\CC$ spanned by the elements $\qlm$ with $\mu=\phi$ can be viewed as the algebra of symmetric functions in variables $x_1,\ldots,x_N$, for large $N$. The elements $h_n$ play the role of the complete symmetric functions and then  $Q_{\gl,\phi}$ 
 corresponds to the classical Schur function $s_\lambda$, expressed as a polynomial in $\{h_i\}$ via the Jacobi-Trudy formula. Determinants similar to the general case for $\qlm$ are used by Koike \cite{Koike} in giving universal
formulae for the irreducible characters of rational representations of
$GL(N)$, along with interpretations in terms of skew Schur functions.
\end{remark}

\subsection*{Examples}

The simplest example is  where $n=p=1$, so that $\gl$ and $\mu$ each have one part of length $1$. In this case the formula gives $\qlm=h_1h^*_1-1$, so that the knot invariant $<K*Q_{\lambda,\mu}>$ is very nearly the reverse-parallel invariant in this case. 

For the figure-eight with zero framing  when $|\gl|=|\mu|=1$ we have $a(\lambda,\mu)=3-2z^2-6z^4-2z^6+(v^2+v^{-2})(-2-z^2-2z^4+z^6)+(v^4+v^{-4})(1+2z^2+z^4).$ The matrix of coefficients is displayed  below, along with the  invariant for the trefoil with some choice of framing - change of framing involves simply factors of $v^2$.
\medskip

\begin{tabular}{c c c} \emph{Figure eight invariant} &\qquad & \emph{Trefoil invariant} \\[2mm]
{$\begin{array}{rr|rrrrr}
&v&-4&-2&0&2&4\\
z&&&&&&\\
\hline &&&&&&\\[-1.5mm]
6&&&1&-2&1&\\4&&1&-2&-6&-2&1\\
2&&2&-1&-2&-1&2\\
0&&1&-2&3&-2&1
\end{array}$}
& \qquad &
{$\begin{array}{rr|rrrr}
&v&0&2&4&6\\
z&&&&&\\
\hline &&&&&\\[-1.5mm]
&&&&&\\
4&&&1&-2&1\\2&&1&2&-7&-4\\
0&&&1&-4&4
\end{array}$}
\end{tabular}

\subsection*{Relations with the Kauffman polynomial}
 In \cite{Rudolph} Rudolph demonstrated a relation between the Kauffman polynomial of a link and the Homfly reverse parallel invariant. His exact result can be described by using the decoration $Q_{\gl,\gl}$ with $|\gl|=1$, as above, on all components of a link $L$. Then the Homfly polynomial of this decorated link determines an element of ${\bf Z}_2[v^{\pm1},z^{\pm1}] $ when the coefficients are reduced  mod  $2$. Rudolph showed that this invariant is the same as the Kauffman polynomial of the link, again with coefficients reduced mod  $2$, when the Kauffman variables $v$ and $z$ are replaced by $v^2$ and $z^2$, and both Kauffman and Homfly are normalised to have the value $1$ on the empty diagram.  The $(1,1)$-tangle invariants above should then reduce to the Kauffman polynomials of the figure eight or trefoil knots, normalised to have the value $1$ on the unknot, with this change of variable. It is reassuring to compare the mod  $2$ reduction of the invariants above with the coefficients for the Kauffman polynomials of these knots shown below, \cite{knotinfo}.

\medskip
 
\begin{tabular}{c c c} \emph{Kauffman polynomial for figure eight} &\qquad & \emph{Kauffman polynomial for trefoil}\\[2mm]
{$\begin{array}{rr|rrrrr}
&v&-2&-1&0&1&2\\
z&&&&&&\\
\hline &&&&&&\\[-1.5mm]
3&&&1&&1&\\2&&1&&2&&1\\
1&&&-1&&-1&\\
0&&-1&&-1&&-1
\end{array}$}
& \qquad &
{$\begin{array}{rr|rrrr}
&v&-5&-4&-3&-2\\
z&&&&&\\
\hline &&&&&\\[-1.5mm]
&&&&&\\
2&&&1&&1\\1&&1&&1&\\
0&&&1&&2
\end{array}$}
\end{tabular}

\subsection*{A possible extension}

 Blanchet and Beliakova  \cite{Blanchet} describe  a decoration $y_\gl$ in the Kauffman skein of the annulus corresponding to each partition $\gl$. Together these account for all possible Kauffman satellite invariants. Where an unoriented link is decorated by one such element $y_{\gl_i}$ on each component its Kauffman polynomial may be compared with the Homfly polynomial of the same link decorated correspondingly by the elements $Q_{\gl_i,\gl_i}$.   The invariant for decorations $y_\gl$ and $Q_{\gl,\gl}$ requires the use of the parameter $s$ with $z=s-s^{-1} $  unless the partition $\gl$ is self-conjugate. When working mod $2$, replacing $s$ by $s^2$ will also have the effect of replacing $z$ by $z^2$.
 Limited evidence suggests the following extension of Rudolph's result from the case $|\gl|=1$  to general Kauffman satellite invariants.
\begin{conjecture} Decorate each component $L_i$ of a framed
 unoriented link $L$ by $y_{\gl_i}$.  The Kauffman polynomial of this decorated link, with $v,s$ replaced by $v^2,s^2$ and the coefficients reduced mod $2$, equals the mod $2$ reduction of the Homfly polynomial of $L$ when each $L_i$ is decorated by $Q_{\gl_i,\gl_i}$.
\end{conjecture}

Known results about quantum dimensions allow the conjecture to be confirmed for the unknot, and for the meridian maps. It is possible that this information can be combined with the branching rules for multiplying $y_\gl$ and $Q_{\gl,\gl}$ by single strings in their respective skeins to give a proof of the conjecture. It would certainly be of interest to study further  the $(1,1)$-tangle invariants for $Q_{\gl,\gl}$.

\leftline{Department of Mathematical Sciences}
\leftline{University of Liverpool}

\end{document}